\documentclass[a4paper,11pt]{elsarticle}
\usepackage[left=2cm,top=2cm,bottom=2cm,right=2cm]{geometry}
\usepackage{amssymb}
\usepackage{amsmath}
\usepackage{amsthm}
\usepackage{epsfig}
\usepackage{float}

\usepackage{multicol}
\usepackage{enumerate}

\theoremstyle{definition}
\newcounter{dummy} \numberwithin{dummy}{section}
\newtheorem{theorem}[dummy]{Theorem}

\begin{document}

\begin{frontmatter}

\title{{\Large \bf Nonlinear oscillators via \v{C}eby\v{s}\"ev quintic approximations}}

\author[MB]{Martina Boschi}

\author[DR]{Daniele Ritelli\corref{COR}}

\author[GS]{Giulia Spaletta}

\cortext[COR]{Corresponding author: Daniele Ritelli,
Department of Statistical Sciences, University of Bologna, 
Via Belle Arti 41,
40126 Bologna, Italy.
Email: daniele.ritelli@unibo.it}

\address[MB]{\em
University of Italian Switzerland,
Faculty of Informatics.
martina.boschi@usi.ch
ORCID 0000--0003--3473--729X}

\address[DR]{
Department of Statistical Sciences, University of Bologna.
daniele.ritelli@unibo.it
ORCID 0000--0001--8805--8132}

\address[GS]{
Department of Statistical Sciences, University of Bologna.	
Member of INdAM--GNCS.
giulia.spaletta@unibo.it
\\
	ORCID 0000--0002--6871--0864}



\begin{abstract}
Aim of this work is the study of differential equations
governing non--dissipative non--linear oscillators;
these arise in different physical models
such as the treatment of relativistic oscillators, 
from the first contribution due to 
\cite{maccoll1957}
and the further analysis in
\cite{mickens1998},
up to generalizations to Duffing's relativistic oscillators
\cite{younesian2012};
they also appear in non--relativistic models
as that in 
\cite{sun2007},
which deals with cables with an attached midpoint mass,
or some harmonic Duffing oscillators discussed in 
\cite{mickens2001}, 
\cite{van2019} 
and 
\cite{razzak2016analytical}. 
From an exquisitely mathematical viewpoint, 
all these models, 
further than describing the one--dimensional motion of a particle,
share being governed by the autonomous equation 
$\,\ddot{x}=f(x)\,$ 
where the restoring force is an odd function, 
and the consequent problem of inverting the associated time--integral; 
the latter can rarely be solved in explicit terms, 
excluding the well--known (both unforced) 
cases of pendulum and Duffing equations.

\noindent
In this paper the inversion issue is treated by 
near--minimax approximation of the restoring force
via fifth--order \v{C}eby\v{s}\"ev polynomials,
on a normalised integration interval: 
this allows time--integral inversion
for the (generalised Duffing) quintic oscillator; 
in fact,
the  particular choice of orthogonal polynomials
turns out to be very appropriate in yielding
an approximate normalised system
whose solutions effectively represent
those of the original problem;
moreover,
when restoring forces are described by odd functions,
the approximate equations are solvable in closed form
via Jacobian elliptic functions. 
\end{abstract}



\begin{keyword}
Non--linear oscillators
\sep
\v{C}eby\v{s}\"ev 
polynomials
\sep
functions approximation
\sep
elliptic integrals
\sep
Jacobian elliptic functions
\sep
\MSC 
33C45
\sep
33E05
\sep
34A05
\sep
65D15
\sep
68W30
\sep
70K75
\end{keyword}

\end{frontmatter}

\section{Introduction}
\label{sec:intro}

\noindent
In this paper 
we present a method for determining analytic solutions 
to fifth--order approximations of non--linear oscillatory systems
governed by odd functions, 
and consequently with even potential energy.
The method described here is based on previous results 
regarding exact analytic solutions of quintic oscillators, 
due to 
\cite{belen2017}, 
\cite{belen2016}, 
\cite{citterio2009}, 
\cite{zuniga2013}, 
\cite{ritelli2010},
which are used by 
\cite{jonckheere1971}, 
and later by 
\cite{zuniga2014} 
and
\cite{biglabo2019},
to treat approximate models
obtained via \v{C}eby\v{s}\"ev polynomial up to degree five.
Our method is applied to the relativistic oscillator 
proposed by McColl 
\cite{maccoll1957} 
and then studied in depth, with different techniques,
in 
\cite{mickens1998}, 
\cite{belen2009}, 
\cite{biazar2014}, 
\cite{hosen2016}. 
We focus on solving non--linear differential equations,
coupled to the originals ruling these models,
using \v{C}eby\v{s}\"ev approximants to their restoring force.

\noindent
Among the techniques
related to the study of non--linear oscillatory phenomena, 
we mention briefly the most popular:
Lindstedt--Poincar{\'e} perturbation methods; 
multiple time--scale methods
\cite{nayfeh1973}, 
\cite{nayfeh1979}, 
\cite{mickens1996};
the generalised averaging method of
Krylov, Bogoliubov and Mitropolski 
\cite{krylov1949}, 
\cite{mickens1996}; 
the approximate variational method,
or energy--balance
\cite{he2002}, 
\cite{belendez2008a},
to evaluate angular frequencies of non--linear oscillators; 
the harmonic--balance method
\cite{mickens1996}, 
\cite{mickens2010}
\cite{gottlieb2004}, 
\cite{gottlieb2006}, 
\cite{wu2006}, 
\cite{belen2008}.  
A notable source for Duffing oscillators is \cite{kovacic2011},
while, for an overview of all these methods,
we highlight 
\cite{cveticanin2014}, 
\cite{mickens2010} 
and again 
\cite{kovacic2011}.  
The problem of period--amplitude dependence was analyzed in 
\cite{rothe1985} 
through classical thermodynamic equilibrium theory, 
and an asymptotic estimate of period is obtained
for the particular case of the
predator--prey Volterra--Lotka model,
which, 
after a suitable change of variable,
is a conservative Hamiltonian system.
This approach is extended to a wide class of Hamiltonian 
non--dissipative system
in 
\cite{ritelli2004},
via Laplace transform and asymptotic expansions. 

\noindent
The differential equations examined in this work are as follows:

\noindent
\begin{minipage}{0.25\linewidth}
\begin{equation}
\ddot{x}
=
-\dfrac{x}{\sqrt{1+x^2}}\;,
\label{iii)}
\end{equation}
\end{minipage}
\begin{minipage}{0.35\linewidth}
\begin{equation}
\qquad
\ddot{x}
=
-x-\dfrac{b\;x}{\sqrt{1+x^2}}\;,
\quad
\label{v)} 
\end{equation}
\end{minipage}
\begin{minipage}{0.35\linewidth}
\begin{equation}
\qquad
\ddot{x}
=
-x-x^3-\dfrac{b \;x}{\sqrt{1+x^2}} \,.
\quad
\label{vi)}
\end{equation}
\end{minipage}
\vskip 9pt

\noindent
Equation
(\ref{iii)}) is related to the relativistic oscillator
introduced in 
\cite{maccoll1957}  
and studied in depth in 
\cite{mickens1998}.
It is actually obtained from 
$\,\ddot{x}+(1-\dot{x}^2)^{3/2}\;x=0\,$
studied in phase--space, 
after a change of variable;
details
are 
well--known and 
reported in several papers, 
such as	the already mentioned
\cite{maccoll1957},  
\cite{mickens1998}, 
and also 
\cite{azami2009},
\cite{belen2009},
\cite{belendez2008higher}.

\vskip 5pt
\noindent
Dynamics of cables 
with an attached midpoint mass
are modeled by (\ref{v)}).
We highlight the contributions of 
\cite{jamshidi2010application}, 
\cite{zhao2009he},
\cite{belendez2007application}, 
\cite{marion2013classical}, 
\cite{younesian2012},
\cite{mickens2001},
\cite{razzak2016analytical}, 
where classical approximate analytic methods
are employed
through some algebraic procedures,
such as an adapted variant of harmonic--balance.

\vskip 5pt
\noindent
The Duffing relativistic oscillator
(\ref{vi)}) 
is treated in
\cite{younesian2012}
using He’s energy--balance method. 

\vskip 5pt
\noindent
Differential equations
(\ref{iii)})--(\ref{vi)}) 
are all of the form 
$\,\ddot{x}=f(x)\,,$
where the restoring force $\,f\,$ is an odd continuous function. 
Assuming 
motion starts from rest, 
i.e. $\,\dot{x}(0)=0\,,$ 
and choosing an initial displacement $\,a>0\,$ 
so that $\,f(x)\neq 0 \quad \forall x\in ]0\,,a]\,,$
the resulting motion is periodic 
and 
the 
particle
satisfies 
$\,-a\leq x(t)\leq a\;,\; \forall\, t\in \mathbb{R}\,.$
In other words,
we study 
an initial value problem (IVP) 
of the form:
\begin{equation}
\begin{cases}
\ddot{x}=f(x)\;,
\\
x(0)=a\;,\quad 
\dot{x}(0)=0\;,
\end{cases}
\label{mickg}
\end{equation}
where $\, f:[-a\, ,a]\to\mathbb{R}\, $ is continuous 
and such that $\, f(-x)=-f(x)\quad \forall  x\in ]0\, ,a]\,;$ 
without loss of generality, 
$\,f(a)<0\,$
can be assumed.
Now, let us introduce the even function:
\begin{equation}
\Phi(x):= - 2\int_x^a f(s)\,{\rm d}s\;.
\label{eqn:phi}
\end{equation}
It is $\,\Phi(\pm a)=0\,,$ 
where 
both zeros are simple 
in the cases of our interest. Moreover:
\begin{equation}
\mathbb{T}= 
2 \displaystyle \int_{-a}^a \; \dfrac{1}{\sqrt{\Phi(s)}}\;\; {{\rm d}s}
\label{main:gensolper}
\end{equation}
is the period of the 
solution to (\ref{mickg}).
This solution
is implicitly defined for 
$\,|x|\leq a\,$
by the
time--integral: 
\begin{equation}
t \;\;=\;\Psi(x)\;,
\qquad
\qquad
\qquad
\Psi(x)\;:=\;\;
\displaystyle
\int_x^a\; \dfrac{1}{\sqrt{\Phi(s)}} \;\; {{\rm d}s}\;.
\label{main:gensol}
\end{equation}

\noindent
At this point,
methods 
of approximation
are necessary
since the
integral in (\ref{main:gensol})
can rarely be evaluated first in closed form and then inverted,
to yield
$\;
x \;\;=\;\Psi^{-1}(t)\;.
\;$

\noindent
Aim of this work
is, namely,
to provide the explicit solution, 
expressed via Jacobian elliptic functions, 
of an approximate problem,
obtained 
by substituting 
the restoring force
$\, f(x)\,$
with
its fifth--order \v{C}eby\v{s}\"ev polynomials
of the first kind,
due to their capability to provide good functions approximation.

\noindent
As it is well--known,
\v{C}eby\v{s}\"ev's 
are a numerable family of polynomials,
orthogonal
with respect to the weigth function 
$\,w(u)=1/\sqrt{1-u^2}\,,$
and defined 
for $\, -1 \leq u \leq 1\,$
by the formulae:
\begin{equation*}
T_n(u)=
\cos\Big(n\;\;\arccos(u)\Big)
=
\;
_{2}{\rm F}_1  
\left(
\begin{matrix}
n \;, -n\\
  \dfrac12
\end{matrix}
\;\left|\;
	\dfrac12\, (1-u) \right.
\right)
\qquad
\qquad
n\in \mathbb{N}\;.
\end{equation*}

\noindent
Here 
$\, _2{\rm F}_1\, $ denotes the Gauss hypergeometric function
\cite{graham-knuth-patashnik}.

\noindent
\v{C}eby\v{s}\"ev
polynomials 
form a complete orthogonal set
on $\,[-1\,,1]\,$ 
in the appropriate Sobolev space, 
thus a function  $\,g(u)\,$
can be expressed
on its domain
$\,[-1\,,1]\,$
via the
expansion:
\begin{equation}
g(u) =
\tilde g(u)
+
E_r(u)\;,
\qquad
\tilde g(u)
:=
\dfrac12\;\alpha_0 + \displaystyle \sum_{n=1}^r\; \alpha_n\;T_n(u)\;,
\qquad
E_r(u)
:=
\displaystyle \sum_{n=r+1}^{\infty}\; \alpha_n\;T_n(u)
\;.
\label{eqn:ChebTruncatedSeries}
\end{equation}

\noindent
If $\, g\,$ is Lipschitz continuous on $\,[-1\,,1]\,,$
then it has a unique
representation 
as the infinite
\v{C}eby\v{s}\"ev
series
(\ref{eqn:ChebTruncatedSeries}),
which is absolutely 
and uniformly 
convergent, with coefficients
defined using the weighted inner product
\cite{trefethen}:
\begin{equation*}
\displaystyle
\alpha_0=
\dfrac{1}{\pi}\;
\displaystyle
\int_{-1}^{1}\;
\dfrac{1}{\sqrt{1-s^2}} \;\; T_n(s)\;\; g(s)\;\;{\rm d}s\;,
\qquad
\quad
\displaystyle
\alpha_n=
\dfrac{2}{\pi}\;
\displaystyle
\int_{-1}^{1}\;
\dfrac{1}{\sqrt{1-s^2}} \;\; T_n(s)\;\; g(s)\;\;{\rm d}s\;\,
\quad 
\mbox{for}
\quad 
n\ge1\,.
\label{eqn:coeff}
\end{equation*}

\noindent
Recall that
it is
$
\, |T_n(u)| \le 1\;\; \forall u\in [-1\;,1]\,. $
Moreover,
$\, T_n(u)\,$ 
has $\,n\,$ distinct real roots 
in $\, ]-1\,,1[\,,$
and
it
has
$\,n+1\,$ extrema
in $\, [-1\,,1]\,$
at which
it
takes 
alternating values $\, \pm 1\,.$
Thus,
if
coefficients 
$\,\alpha_n\,$
decrease in magnitude sufficiently rapidly
(which depends on regularity of $\,g\,$),
then
$\;
E_r(u)
\simeq
\alpha_{r+1}\;T_{r+1}(u) \; $
equioscillates $\, r+2\,$ times on $\,[-1\,,1]\,,$
implying that
$\, \tilde g\,$
is  a near--minimax approximant for $\,g\,$
\cite{phillips-taylor}.
Coefficients 
$\,\alpha_n\,$ 
can be determined explicitly for some functions,
otherwise they
need discretisation via quadrature formulae.
Even so, 
among methods yielding minimax or near--minimax approximations, 
\v{C}eby\v{s}\"ev
series are effective and easy to handle.

\vskip 2pt
\noindent
To apply 
\v{C}eby\v{s}\"ev 
approximation to
the nonlinear oscillators under study,
the displacement $\,a\,$
is normalised
to the interval $\,[-1\;,1]\;$
via a change of dependent variable $\,u=x/a\,,$
and 
the following
equivalent
IVP 
is considered in place of
(\ref{mickg}):
\begin{equation}
\begin{cases}
\ddot{u}= f_a(u) \;,
\qquad\qquad\qquad 
\qquad\qquad\qquad 
 f_a(u) :=
\dfrac{1}{a}\;\; f(a\,u)\;,
\\
u(0)=1\;, \quad
\dot{u}(0)=0\;.
\end{cases}
\label{main:genuni}
\end{equation} 

\vskip 2pt
\noindent
\noindent
We then choose to 
describe
the normalised restoring force
$\, f_a(u)\,,$
which is an odd function,
in terms of polynomials $\,T_n(u)\,;$
for our purposes, 
$\,f_a(u)\,$ 
is expanded
in \v{C}eby\v{s}\"ev series 
truncated (or projected) at fifth--order:
\begin{equation}
f_a(u)\, \simeq\, 
\tilde f_a(u) :=
\alpha_1\, T_1(u)+\alpha_3\, T_3(u)+\alpha_5\, T_5(u)\;,
\label{eqn:approximatefa}
\end{equation}
where: 
\begin{equation*}
T_1(u)=u\;,
\qquad
T_3(u)=-3\,u+4\, u^3\;,
\qquad
T_5(u)=5\,u-20\,u^3+16\, u^5\;,
\label{eqn:T1-T3-T5}
\end{equation*}
and
\begin{equation}
\alpha_n=
\dfrac{2}{\pi}\;
\displaystyle
\int_{-1}^{1}\;
\dfrac{1}{\sqrt{1-s^2}}\;\;T_n(s)\;\; f_a(s)\;\;{\rm d}s\;,
\qquad
\qquad n=1\,,3\,,5\,.
\quad
\label{eqn:weighted}
\end{equation}

\noindent
Expressing the approximate normalised 
force
$\, \tilde f_a(u) \,$
in the monomial base,
a new
IVP replaces
(\ref{main:genuni}):
\begin{equation}
\begin{cases}
\ddot{u}=-(c_1\; u+c_3\;u^3+c_5\;u^5)\;,\\
u(0)=1\;,\quad
\dot{u}(0)=0\;,
\end{cases}
\label{main:appr}
\end{equation}
where, 
setting
$\; \mathcal{C} = -2^5/\pi\,:$
\begin{equation}
\begin{aligned}
c_1 = -( \alpha_1-3\,\alpha_3+5\, \alpha_5)
&=&
\mathcal{C}\;
\displaystyle
\int_{-1}^1\;
\dfrac{1}{\sqrt{1-s^2}}
\quad
&(\dfrac{35}{16}\, s-7\, s^3+5\, s^5)
&f_a(s)
\quad 
&{\rm d}s\;,
\\[3pt]
c_3=-4\,(\alpha_3-5\,\alpha_5)
&=&
\mathcal{C}\;
\displaystyle
\int_{-1}^1\; 
\dfrac{1}{\sqrt{1-s^2}}
\quad
&(-7\, s+26\, s^3-20\, s^5)
&f_a(s)
\quad 
&{\rm d}s\;,
\\[3pt]
c_5=-16\,\alpha_5
&=&
\mathcal{C}\;
\displaystyle
\int_{-1}^1\; 
\dfrac{1}{\sqrt{1-s^2}}
\quad
&(5\, s-20\, s^3+16\, s^5)
&f_a(s)
\quad 
&{\rm d}s\;.
\end{aligned}
\label{koeffy}
\end{equation}

\noindent
In 
$\S$~\ref{sec:general5oscillator},
a novel solution procedure is 
presented for 
the general quintic oscillator
\eqref{main:appr} 
in terms of Jacobian elliptic functions, leading
to the determination of
exact periods and frequencies.
In 
$\S$~\ref{sec:application2relativistic-oscillator},
the solution process is illustrated on 
the relativistic oscillator
(\ref{iii)}):
its normalised and 
quinticated approximation
is built
and the exact integration obtained
in 
$\S$~\ref{sec:general5oscillator}
is applied
to solve it;
quality of the
results
obtained
is also
validated.
In
$\S$~\ref{sec:5oscillators},
a similar application to 
oscillators (\ref{v)}) and (\ref{vi)})
proves both feasibility and robustness of the solution process introduced.
The conclusive
$\S$~\ref{sec:conclusion} reports
some final comments and indications for future work.

\section{General quintic oscillator}
\label{sec:general5oscillator}

\noindent
Consider the family of 
IVPs~(\ref{main:appr}).
We highlight several contributions 
for this kind of problems 
due to 
\cite{belen2016}, 
\cite{belen2017}, 
\cite{belen2012}, 
\cite{belen2016b}, 
\cite{citterio2009},
\cite{khalil2021},
\cite{ritelli2010},
\cite{zuniga2013}. 
Application of 
\eqref{eqn:phi}
and
\eqref{main:gensol}
to 
the approximate normalised force
in \eqref{main:appr}
shows that
the (squared)
solution of this IVP
is based on the evaluation of an elliptic integral:
\begin{equation}
t=
\sqrt{\frac{3}{2}}\;\;
\displaystyle
\int_{u^2}^{1}\;\;
\dfrac{1}{\sqrt{s\; \left(1-s\right) \;h_2(s)}}\quad
{{\rm d}s}
\;,
\label{quinyg}
\end{equation}
with
\begin{equation*}
h_2(s)=(6\, c_1+3\,c_3+2\,c_5)+\left(3\, c_3+2\,c_5\right)\,\,s+2\, c_5\, s^2\;.
\label{eqn:h2}
\end{equation*}
The discriminant of polynomial $\, h_2(s)\, $ is,
discarding a factor of value $\,3\,:$
\begin{equation}
\Delta=3\, c_3^2-4\, c_5\, \left(4\, c_1+c_3+c_5\right)\;.
\label{discr}
\end{equation} 

\noindent
Given the physical nature of the restoring forces 
acting in the models of interest,
coefficients 
$\,c_1\,, c_3\,,c_5\,$ 
can be assumed to be
such that  
$\, h_2(s) > 0 \quad \forall s\in(0\,,1)\,.$ 
This property is assured if 
$\,c_5>0\,$ 
together with one of the two conditions:
\begin{equation}
\mbox{(i)}
\quad \Delta \leq 0\;; 
\qquad
\qquad
\qquad
\mbox{(ii)}
\quad
\Delta > 0\,
\quad
\mbox{and}
\quad
 \,6\, c_1+3\,c_3+2\,c_5>0\,. 
\qquad
\label{eqn:conditions} 
\end{equation}

\begin{figure}[h]
\begin{center}
\includegraphics[width=220pt]{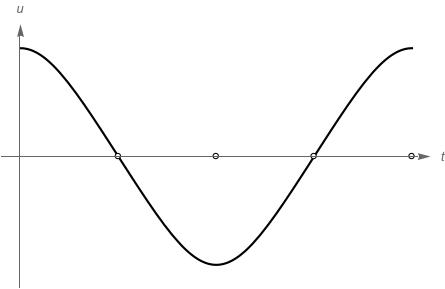}
\caption{{Plot of the solution to IVP
\eqref{main:appr}
for $\,c_1=1\;,c_3=2\;,c_5=3\,,\,$ in one period interval}}
\label{fig:0}
\end{center}
\end{figure}

\noindent
The following
Theorems 
\ref{ddue} and \ref{unnno}
provide closed--form solution and period
for \eqref{main:appr},
respectively
under
conditions (i) or (ii). 
Notice that, in the latter case, 
the 
roots of $\,h_2(s)\,,$ 
further than being 
real and distinct,
are
both negative, 
due to Descartes' rule of signs.
Solution of
\eqref{main:appr}
has a cosine wave behaviour,
as
Figure
\ref{fig:0}
illustrates for the example case 
$\,c_1=1\;,$ $ c_3=2\;,$ $ c_5=3\;.$

\begin{theorem}
\label{ddue}
Given
the time--integral
\eqref{quinyg},
assume
$\,c_5>0\,$ 
and 
$\,\Delta \leq 0\,$ 
in \eqref{discr},
and define:
\begin{equation}
A =
\dfrac{\sqrt[4]{6}}{2} \; \;
\dfrac{1}{\sqrt[4]{\mathcal{P} \;\; \mathcal{Q}}} \;,
\qquad
\qquad
B = 
\dfrac{1}{6}\;
\dfrac{\mathcal{Q}}{\mathcal{P}}\;,
\qquad
\qquad
\qquad
k^2 = 
\dfrac12-
\dfrac{\sqrt{6}}{8}
\;
\dfrac{ \mathcal{K}} {\sqrt{\mathcal{P}\;\; \mathcal{Q}}}\;,
\label{parram}
\end{equation}
with
\begin{equation}
\mathcal{P}
=
c_1+c_3+c_5\;,
\qquad
\qquad
\mathcal{Q}
=
6\, c_1+3\, c_3+2\, c_5\;,
\qquad
\qquad
\mathcal{K}
=
4 \,c_1+3 \,c_3+2 \,c_5\;.
\label{eqn:parametri}
\end{equation}

\noindent
Then, the solution 
of 
IVP~\eqref{main:appr}
is:
\begin{equation}
u^2(t)=
\dfrac{\sqrt{B}}
{
\sqrt{B}+
\cot^2
\Bigg(
\dfrac{1}{2}\;\; 
{\rm am}\left(2\; \mathbf{K}(k) -\dfrac{t}{A}\;,\;k\right)
\Bigg)
}\;.
\label{udue}
\end{equation}
where
$\,{\rm am}(s\,,k)\,$ 
indicates the Jacobi amplitude function,
i.e.
the inverse of the elliptic 
integral of first kind 
$\, F(\varphi\,,k)\,,$
meaning that
$\, \varphi={\rm am}(s\,,k) \;$
iff
$\;
s
=
F(\varphi\;,k)\;, $
where:
\begin{equation}
F(\varphi\;,k) \;:=\;
\displaystyle
\int_0^{\sin(\varphi)}\;
\dfrac{1}
{
\sqrt{
(1-s^2)\;(1-k^2\;s^2)
}
}
\;
\;
{\rm d}s\;,
\qquad
\qquad
-\dfrac{\pi}{2}<\varphi<\dfrac{\pi}{2}\;,
\label{eqn:EllipticIntegralF1}
\end{equation}
while
$\, \mathbf{K}(k)\, \;:=\; F(\dfrac{\pi}{2}\;,k)\; $ 
 denotes 
the complete elliptic integral of first kind, 
with elliptic modulus $\,k\,.$ 

\vskip 2pt
\noindent
Solution 
\eqref{udue}
is periodic, with
period:
\begin{equation}
\mathbb{T}=8\,\;A\,\;\mathbf{K}(k)\;,
\label{main:appr:per2}
\end{equation}
and it is positive 
for $\,0 \leq t \leq \dfrac{1}{4}\,\mathbb{T}\;,\quad
\dfrac{3}{4}\, \mathbb{T}\leq t\leq \mathbb{T}\,,$ 
while it is negative 
for
$\, \dfrac{1}{4}\, \mathbb{T}< t<\dfrac{3}{4}\,\mathbb{T}\,.$
\end{theorem}

\vskip 2pt
\begin{proof}
The integral in \eqref{quinyg} 
can be evaluated 
using
entry
3.145--2
in
\cite{grad2000}, 
recalled here:
\begin{equation}
\int_\beta^v
\;
\dfrac{1}
{\sqrt{(\eta-s)\;(s-\beta)\;\;\big((s-m)^2+n^2\big)}}
\;
{\rm d}s\;
=
\;
\dfrac{1}{\sqrt{p\, q}}\,F(\varphi(v)\;,k)\;,
\label{259:00}
\end{equation}
where: 
\begin{equation*}
\beta<v<\eta\;,
\qquad
\quad
p^2:=(m-\eta)^2+n^2\;, 
\quad
\qquad q^2:=(m-\beta)^2+n^2 \,,
\label{eqn:259bis}
\end{equation*}
and
\begin{equation*}
\label{259:00a}
\displaystyle
{\varphi(v)=2\,{\rm arccot}\sqrt{\dfrac{q\;(\eta-v)}{p\;(v-\beta)}}}\;,
\qquad
\displaystyle
{k^2=\dfrac14\,\dfrac{(\eta-\beta)^2-(p-q)^2}{p\;q}}\;.
\end{equation*}
	
\vskip 2pt
\noindent
In the case of \eqref{quinyg}, 
$\,\beta=0\;,\,\eta=1\;, v=u^2\;, $ 
and polynomial $\, h_2(s)\,$
is rearranged as:
\begin{equation*}
\dfrac{1}{2\;c_5}\;\; h_2(s) \;=\;
\Bigg(
s
\;+\; \Big(\dfrac{3\,c_3}{ 4\;c_5 } \;+\;\dfrac12 \Big) \Bigg)^2
\;+\;
\Bigg(
\dfrac{ 3\, c_1}{ c_5 }
\;-\;
\Big(\dfrac{ 3\, c_3}{4\; c_5 }\Big)^2
\;+\;
\dfrac{ 3\, c_3}{ 4\; c_5 }
\;+\;
\dfrac34
\Bigg)\,.
\label{h2-rearranged}
\end{equation*}

\vskip 2pt
\noindent
To apply 
formula \eqref{259:00}, 
the
integral in \eqref{quinyg} 
must further be rewritten
as the difference of integrals of the same integrand 
on intervals $\,[0\,,1]\,$
and $\,[0\,,u^2]\,.$
Equation \eqref{quinyg} thus becomes:
\begin{equation}
t=
2\; \;A\;\; \mathbf{K}(k)
- A\;\;
F\Bigg(2\; {\rm arccot}\left(\sqrt[4]{B}\;\; \sqrt{\Upsilon(u^2)}\right)\;,\; k\Bigg)\;,
\label{eqn:time1}
\end{equation}
where 
$\,\Upsilon(v)=(1-v)/v\,.$
At this point, 
due to invertibility of the elliptic integral of first kind,
inversion of the time--integral equation
\eqref{eqn:time1}
is possible,
and it yields
solution
\eqref{udue}.
In a similar way, 
\eqref{main:appr:per2} can be proved to provide
the motion period.  
\end{proof}

\vskip 5pt
\begin{theorem}
\label{unnno}
Given
the time--integral
\eqref{quinyg},
assume
$\,c_5>0\,$ 
and 
condition (ii) in
\eqref{eqn:conditions},
that is
$\; h_2(s)=(s-s_1)\,(s-s_2)\;$
where
$\, s_1<s_2<0\,$ 
are its
real roots.  
Now, define:
\begin{equation}
k^2=\dfrac{s_2-s_1}{s_1 \left(s_2-1\right)}\;.
\label{main:alfakappa}
\end{equation}

\noindent
Then, the solution of IVP~\eqref{main:appr} is:
\begin{equation}
u^2(t)=
s_1+
\dfrac{s_1\;(s_1-1)}
{{\rm sn}^2\left(\sqrt{\dfrac{c_5\;s_1\; 
\left(s_2-1\right)}{3}}\;\;\;t\;,\;\;k\right)
-s_1}\;,
\label{uuno}
\end{equation}
where
$\, {\rm sn}(s\,,k) \,$
is the Jacobi sine amplitude function, i.e.
$\, {\rm sn}(s\,,k) = \sin(\varphi) \,,$ with $\, \varphi={\rm am}(s\,,k)\,.$

\vskip 2pt
\noindent
Solution~\eqref{uuno} is periodic, with period:
\begin{equation}
\mathbb{T}=
\dfrac{4\;\sqrt{3}}{\sqrt{c_5\;s_1\;\left(s_2-1\right)}}
\;\; \mathbf{K}\left(k\right)\;,
\label{main:appr:per1}
\end{equation}

\noindent
and 
it is positive 
for $\,0 \leq t \leq \dfrac{1}{4}\,\mathbb{T}\;,\;\;
\dfrac{3}{4}\, \mathbb{T}\leq t\leq \mathbb{T}\,,$ 
while it is 
negative 
for
$\, \dfrac{1}{4}\, \mathbb{T}< t<\dfrac{3}{4}\,\mathbb{T}\,.$
\end{theorem}

\vskip 2pt
\begin{proof}
To evaluate the integral in \eqref{quinyg},
entry 3.147--7 of \cite{grad2000} is used, 
recalled below:
\begin{equation*}
\int_v^{\eta}
\dfrac{{\rm d}s}
{\sqrt{(\eta-s)\, (s-\beta)\, (s-\gamma)\, (s-\delta)}}
=
\dfrac{2}{\sqrt{(\eta-\gamma)\,(\beta-\delta)}}\;\;F(\varphi(v),k)\;,
\label{eqn:3147-7-GR}
\end{equation*}
where:
\begin{equation*}
\delta<\gamma<\beta\le v < \eta\;,
\qquad
\varphi(v)=
\arcsin
\sqrt{\dfrac{(\beta-\delta)(\eta-v)}{(\eta-\beta)(v-\delta)}}\,,
\qquad
k^2=\dfrac{(\eta-\beta)(\gamma-\delta)}{(\eta-\gamma)(\beta-\delta)}\;.
\end{equation*}
In the case of \eqref{quinyg},
$\,\beta=0\,, \eta=1\,,\, v=u^2\,,\,\gamma=s_2\,,\,\delta=s_1\;.$ 
Thus, the motion period is given by \eqref{main:appr:per1}, while
\eqref{uuno} provides the solution, after the relevant computations,
not reported here, as they are similar to those performed 
in  proving Theorem 
\ref{ddue}.
\end{proof}

\vskip 5pt
\noindent
If $\,\Delta=0\,,$ integral \eqref{quinyg} 
degenerates into an elliptic integral of third kind,
which is tabulated as 
entry 3.138--6 of
\cite{grad2000}. 
Here, the related computation are omitted for two reasons.
First of all,
condition $\,\Delta=0\,$ is linked to a very particular value 
of the initial displacement $\,a\,.$ 
Secondly,
even though 
the appareance of elliptic integrals of third kind
makes it impossible to invert the time--integral
and compute the solution explicitly,
thanks to the continuous dependence on data, 
the relevant solution can be approximated 
at arbitrary precision 
with solutions obtained in Theorems \ref{ddue} and \ref{unnno}.

\section{Application to the relativistic oscillator}
\label{sec:application2relativistic-oscillator}

\noindent
In the case of the relativistic oscillator
ruled by \eqref{iii)} 
the normalised equation is:
\begin{equation}
\begin{cases}
\ddot{u}=-\dfrac{u}{\sqrt{1+a^2 \;\; u^2}}\;,
\\
u(0)=1\,,\quad \dot{u}(0)=0\;.
\end{cases}
\label{ii)n}
\end{equation}

\noindent
Here, 
function
$\, \Phi(u) \,, $ 
defined in   
\eqref{eqn:phi},
becomes:
\begin{equation}
\Phi(u)= 
\dfrac{2}{a^2}\, \Big(\sqrt{1+a^2} - \sqrt{1+a^2\; u^2}\Big)\,,
\label{phi-psi}
\end{equation}
while, after some algebraic adjustments, function
$\, \Psi(u) \,, $ 
defined in 
\eqref{main:gensol},
is:
\begin{equation}
\begin{split}
\Psi(u)
&=
\frac{a}{ \sqrt{2}\;\; \sqrt[4]{1+a^2}}
\;\;
\int_u^1\;
\dfrac{1}{ \sqrt{1-\sqrt{
\dfrac{1+a^2 \; s^2}{1+a^2}}}}\;\; \;{\rm d}s
\\
&=
\dfrac{\sqrt[4]{1+a^2}}{\sqrt{2}}
\;\;
\int_{_{\sqrt{\frac{1+a^2\; u^2}{1+a^2}}}}^1
\qquad
\dfrac{z}{\sqrt{(1-z) \left(z^2-\dfrac{1}{1+a^2}\right)}}\;\;{\rm d}z\;,
\end{split}
\label{temporel1}
\end{equation}
where it is, obviously,
$\; 0 < \frac{1}{\sqrt{1+a^2}} < \sqrt{\frac{1+a^2\; u^2}{1+a^2}} <1\;. $ 
The integral 
in \eqref{temporel1}
can be expressed in explicit form,
for example
via
entry 3.132--5 of
\cite{grad2000}, 
through which
the time equation \eqref{main:gensol} becomes:
\begin{equation}
t
= \sqrt{2}\;\;
\Bigg(
\mathcal{A}\;\;
E\left(\varphi,k\right)-
\dfrac{
	F\left(\varphi,k\right)
	}{\mathcal{A}}
\Bigg)\;,
\qquad
\qquad
\mathcal{A}=
\sqrt{\sqrt{1+a^2}+1}\;,
\label{temporel2}
\end{equation}
being 
$\, F\left(\varphi,k\right)\,$ 
and $\, E\left(\varphi,k\right)\, $ 
elliptic integrals of first and second kind, with:
\begin{equation}
\varphi=\arcsin\sqrt{\frac{\sqrt{1+a^2}-\sqrt{1+a^2\; u^2}}{\sqrt{1+a^2}-1}}\,,
\qquad 
\qquad 
k= \dfrac{\sqrt{1+a^2}-1}{a}\;.
\label{phikappa}
\end{equation}
Note that, for any $\, a>0\;,$ 
the elliptic modulus $\,k\,$ satisfies the requirement 
$\, k<1\;.$ 

\noindent
Calculation of the integral in \eqref{temporel1} further leads to 
the exact determination of the period of oscillation:
\begin{equation}
\mathbf{T}
= 4\, \sqrt{2}\;\;
\Bigg(
\mathcal{A}\;\;
\mathbf{E}(k)-
	\dfrac{
\mathbf{K}(k)	
}{\mathcal{A}}
\Bigg)\;,
\qquad
\qquad
\mathcal{A}=
\sqrt{\sqrt{1+a^2}+1}\;,
\label{perrel}
\end{equation}

\noindent
where $\mathbf{K}(k)$ and $\mathbf{E}(k)$
are the complete elliptic integrals of first and second kind,
respectively,
with modulus $k$ as in \eqref{phikappa}.
The explicit formula 
\eqref{perrel}
for the period 
of the relativistic oscillator is useful in itself, 
and also because it allows a comparison with the period of 
the approximated quintic oscillator we are to obtain.

\noindent
Deriving solution $\,u\,$ from  
\eqref{temporel2}, in fact,
poses
the computational problem 
represented by 
inversion of the time--integral.
We then, instead,
approximate the normalised restoring force in \eqref{ii)n},
using \v{C}eby\v{s}\"ev polynomials and expressing it in the monomial base
through coefficients 
$\, c_1\,, c_3\,,c_5\,$
given by
\eqref{koeffy}:
according to
the sign of the 
$\, \Delta\,$
discriminant 
built on such coefficients,
the seeked solution 
is
\eqref{udue}
or
\eqref{uuno}.

\vskip 5pt
\noindent
Application of formulae 
\eqref{koeffy}
to the considered problem \eqref{ii)n} suggests 
to introduce three elliptic integrals:
\begin{equation*}
J_{n}(a)=
\int_{-1}^1\;
\dfrac{s^n}
{\sqrt{\left(1-s^2\right) \left(1+a^2\; s^2\right)}}\,{\rm d}s\;,
\qquad n=2\,,4\,,6\;.
\end{equation*}
since,
setting
$\; \mathcal{C} =
 2^5/\pi\,:$
\begin{equation}
\begin{aligned}
c_1&= 
\mathcal{C}\;
\Big(\dfrac{35}{16}\, J_2(a)-7\, J_4(a)+5\, J_6(a)\Big)\;,
\\
c_3&=
\mathcal{C}\;
\left( -7\, J_2(a)+26\, J_4(a)-20\, J_6(a)\right)\;,
\\
c_5&=
\mathcal{C}\;
\left(5\,  J_2(a)-20\, J_4(a)+16\, J_6(a)\right)\;.
\end{aligned}
\label{cocoe}
\end{equation}

\noindent
Using entries 236.16 and 331.01--03
of \cite{byrd1971},
it follows:
\begin{equation}
\begin{aligned}
J_2(a)
&=
\dfrac{2}{a^2} \;
\Big( 
\mathcal{J}\;\;\mathbf{E}(h^2) 
\; -\; 
\dfrac{1}{\mathcal{J}}\;\; \mathbf{K}(h^2)
\Big)\;,
\\
J_4(a)
&=
\dfrac{2}{3\; a^4} \; 
\Big(
2\; (a^2-1)\; \; 
\mathcal{J}\;\; \mathbf{E}(h^2) 
\; - \;
(a^2-2)\;\; 
\dfrac{1}{\mathcal{J}}\;\; \mathbf{K}(h^2)
\Big)\;,
\\
J_6(a)
&=
\dfrac{2}{15\, a^6}\;
\Big( 
(8\,a^4-7\,a^2+8)\;\;
\mathcal{J}\; \;\mathbf{E}(h^2)
\;-\;
(4\,a^4-3a^2+8)\;\; 
\dfrac{1}{\mathcal{J}}\;\; \mathbf{K}(h^2)
\Big)\;,
\end{aligned}
\label{geii}
\end{equation}
where 
$\, \mathcal{J}=\sqrt{1+a^2}\,$
and
the elliptic modulus is given by 
$\, h=a/\mathcal{J}\,.$

\noindent
It is thus possible to identify closed--form expressions 
for coefficients 
$\, c_1\,, c_3\,, c_5\, $
of the approximate quintic oscillator,
inserting values 
\eqref{geii}
of integrals $\, J_2\,,J_4\,,J_6\,$ into
\eqref{cocoe}.

\begin{figure}[h]
\begin{center}
\includegraphics[width=220pt]{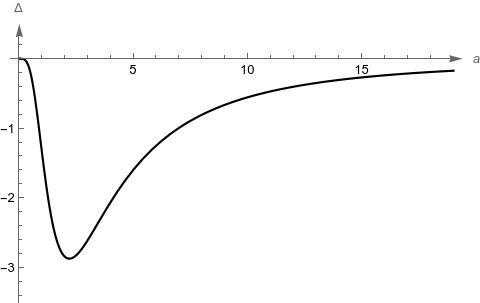}
\caption{Discriminant 
$\, \Delta=\Delta(a)\,$
associated to 
the quinticated
form of 
oscillator
(\ref{iii)}) }
\label{deeltarell}
\end{center} 
\end{figure}

\noindent
The consequent expression
\eqref{discr}
of 
the 
$\, \Delta\,$
discriminant
is complicated,
but still tractable using computer algebra systems, such  as {\sl Mathematica}
\cite{wri},
within which various
graphical tools are also provided,
enabling a visual analysis as that shown in Figure \ref{deeltarell}.

\noindent
At this point,
knowing that $\,\Delta\le0\,,$ 
solution 
and period are those given in 
Theorem \ref{ddue};
condition $\, c_5 >0\, $ can, indeed, also be checked, for example
via built--in visualisation resources.
This whole process indeed involves 
complicated
expressions
which, however, remain 
manageable using computer algebra.

\vskip 5pt
\noindent
What is important to note is that, in order to use the procedure presented, 
what is ultimately needed 
is to calculate the integrals expressing the orthogonal projection
onto the space of \v{C}eby\v{s}\"ev polynomials. 
From that point on, 
by means of the coefficients of the quintic polynomial 
that approximates the normalised restoring force
and by Theorem~\ref{ddue},
or if appropriate 
Theorem~\ref{unnno}, 
one arrives at the approximate solution and its period.

\vskip 5pt
\noindent
To validate the quality of the approximation obtained, 
the computational capacity of {\sl Mathematica}
can again be exploited 
in evaluating the differential operator:
\begin{equation}
L\, u=\ddot{u} \, -\, f_a(u)\,,
\label{operator}
\end{equation}

\noindent
where 
here it is 
$\,  f_a(u) = - 
{u}/{\sqrt{1+a^2\; u^2}}\,,$
since we are studying oscillator (\ref{iii)}),
while 
$\,  f_a(u) = f_a(u\,,b) =
- u -
{b\, u}/{\sqrt{1+a^2\; u^2}}\,$
for
oscillator
(\ref{v)}) 
and
$\,  f_a(u) = f_a(u\,,b)= 
- u - a^2 \, u^3
- {b\, u}/{\sqrt{1+a^2\; u^2}}\,$
for
oscillator 
(\ref{vi)}), 
as we will see in $\S$~\ref{sec-oscillator2} and 
$\S$~\ref{sec-oscillator3}
respectively.

\begin{figure}[h]
\begin{center}
\includegraphics[width=230pt]{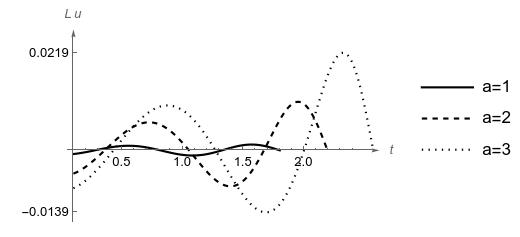}
\includegraphics[width=230pt]{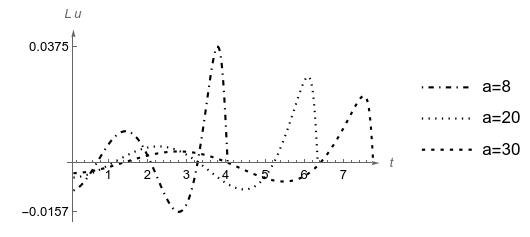}
\caption{Behaviour of $\, L\, u\,,$ 
where $\, u\,$ is solution 
\eqref{udue} with coefficients 
\eqref{cocoe},
for the quinticated 
form		
of
oscillator 
(\ref{iii)})
}
\label{fig-norm123}
\end{center}      
\end{figure}

\noindent
For the quinticated 
IVP associated to
oscillator 
(\ref{iii)}),
Figure
\ref{fig-norm123}
reports
the graph of the deviation from zero  
produced 
in \eqref{operator}
by 
solution
$\, u\,$ 
in the first quarter of the period,
and
at initial displacements equal to
$\, a=1,\,2,\, 3\,$
and
$\, a=8,\,20,\, 30\,;$
plots are kept separate for effective rendering reasons.

\noindent
It
is reasonable
for
the approximation to initially get worse
as the displacement increases, 
given the normalisation of the integration interval.
However,
computation of the maximum difference 
in the first quarter of the period 
yields results in Table
\ref{tab:norm}, showing that an upper bound is
given by 
the value 
$\, 0.0375439\,$
(approximatively and  working in machine precision)
reached
at $\, a=8\,,$ 
from 
where a monotonic decrease can be observed.
\begin{table}[h]
\centering
\begin{tabular}{ccccccc}
\hline
{a }
& 
{1 }
& 
{2 }
&
{3 }
& 
{8 }
& 
{20 }
& 
{30 }
\\
\hline
{$\| L\, u\|_{\infty}$ }
& $\;$ 
{ .0013005 }
$\;$
& $\;$ 
{ .0109030 }
$\;$
& $\;$
{ .0219219 }
$\;$
& $\;$
{ .0375439 }
$\;$
& $\;$ 
{ .0278857 }
$\;$
& $\;$ 
{ .0216839 }
$\;$
\\
\hline
\end{tabular}
\caption{Uniform norm of $\, L\, u\,,$ 
where $\, u\,$ is solution 
\eqref{udue} with coefficients 
\eqref{cocoe}, 
for the quinticated 
form		
of
oscillator 
(\ref{iii)})
}
\label{tab:norm}
\end{table}

\begin{figure}[h]
\begin{center}
\includegraphics[width=220pt]{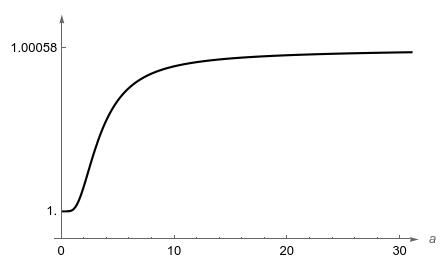}
\caption{Period ratio 
in the case of 
the 
quinticated form of 
oscillator
(\ref{iii)}) }
\label{fig-ratio}
\end{center}       
\end{figure}

\vskip 5pt
\noindent
The validity of the quintic approximations 
is further testified by the ratio 
between exact period 
\eqref{perrel}
of the gravitational oscillator 
and period
\eqref{main:appr:per2} 
of the approximation: 
this ratio
remains close to $\, 1\,,$ 
and bounded above by the value $1.00058\,$ 
for large values of $\,a\,.$
Figure \ref{fig-ratio}
illustrates the ratio behaviour 
for 
displacements 
up to $\,a=30\,.$

\section{Quinticated oscillators}
\label{sec:5oscillators}

\noindent
Here, the procedure introduced in 
$\S$~\ref{sec:general5oscillator} is applied to the conservative
non--linear oscillatory system \eqref{v)} and the Duffing relativistic
oscillator
\eqref{vi)}.
Results obtained are wholly analogous to those seen for the relativistic oscillator
\eqref{iii)}.
In particular, it is still possible to determine closed--form expressions for the 
coefficients of the
\v{C}eby\v{s}\"ev 
quintic approximant, in terms of 
complete elliptic integrals of first and second kind.

\subsection{Nonlinear oscillator (\ref{v)}) }
\label{sec-oscillator2}

\noindent
The normalised IVP,
here, is:
\begin{equation}
\begin{cases}
\ddot{u}=-u-\dfrac{b\;u}{\sqrt{1+a^2 \;\; u^2}}\;,
\\
u(0)=1\,,\quad \dot{u}(0)=0\;,
\end{cases}
\label{eqn-oscillator2normalised}
\end{equation}
thus 
function
$\, \Phi(u) \,, $ 
defined in   
\eqref{eqn:phi},
becomes:
\begin{equation}
\Phi(u)= 1-u^2
+
\dfrac{2\,b}{a^2}\, \Big(\sqrt{1+a^2} - \sqrt{1+a^2\; u^2}\Big)\,,
\label{phi-psi-oscillatore2}
\end{equation}
so that forming
$\, \Psi(u)\,$ as in \eqref{main:gensol}
requires 
$\, 2\, b > \sqrt{1+a^2}-\sqrt{1+a^2\,u^2} > 0\,,$
being $\, a>0\, $ and $\, 0<u<1\,.$

\noindent
The exact period in this case,
obtained
using
entry
3.148--7
in
\cite{grad2000}, 
is:
\begin{equation}
\mathbf{T}
= \dfrac{4}{\sqrt{\mathcal{A}}} \;\;
\Bigg((1+\sqrt{1+a^2})\;\;\;
\mathbf{\Pi}\left(\mathcal{N}\;,\;
\mathcal{K}\right) 
\; - \;\;
\mathbf{K}(\mathcal{K}) \Bigg)
\label{periodo-oscillatore2}
\end{equation}
with
\begin{equation*}
{\mathcal{A}} = \sqrt{1+a^2} + b\;,
\qquad \qquad
\mathcal{N} = \dfrac{1-\sqrt{1+a^2}}{2}\;,
\qquad\qquad
{\mathcal{K}} = \dfrac{\mathcal{N} \;\;(b-\mathcal{N})}{\mathcal{A}}\;.
\end{equation*}

\noindent
$\mathbf{K}(\mathcal{K})$ and $\mathbf{\Pi}(
\mathcal{N}\;,\;\mathcal{K})$
are the complete elliptic integrals of first and third kind,
respectively,
where
$\, \mathcal{N} \,$
is the elliptic characteristic,
while the modulus 
satisfies 
$\, \mathcal{K}<1 \,$ 
for $\, a\,,b>0\,.$

\noindent
For 
system 
\eqref{eqn-oscillator2normalised},
given the exact
computation of integrals 
\eqref{geii}, 
the three coefficients of the
quinticated
approximant 
can also be computed explicitly,
using again 
\eqref{koeffy}, recalled here for reading convenience:
\begin{equation*}
c_1 = -( \alpha_1-3\,\alpha_3+5\, \alpha_5)\;,\qquad
c_3 = -4 \; ( \alpha_3-5\, \alpha_5)\;,\qquad
c_5 = -16 \;  \alpha_5\;,
\end{equation*}
where, setting $\, \mathcal{C}= -2/\pi\;:$ 
\begin{equation}
\begin{aligned}
\alpha_1
&=
\phantom{111}
\mathcal{C}\;   \Big(w_2 + b\; J_2(a)\Big)\;, 
\\
\alpha_3
&= 
-3\; 
\mathcal{C}\; 
\Big(w_2 + b\; J_2(a)\Big) \; 
+ 4 \; 
\mathcal{C}\; 
\Big(w_4 + b\; J_4(a)\Big)\;, 
\\
\alpha_5
&= 
\phantom{-}
5 \;\mathcal{C} \; \Big(w_2 + b\; J_2(a)\Big)
-20 \;
\mathcal{C}\; 
\Big(w_4 + b\; J_4(a)\Big) 
+16 \; 
\mathcal{C}\; 
\Big(w_6 + b\; J_6(a)\Big)\;, 
\end{aligned}
\label{c1-3-5-oscillator1.2}
\end{equation}
with:
\begin{equation}
w_n=  \int_{-1}^1 \;  \dfrac{s^n}{\sqrt{1-s^2}} \;, \qquad n=2\,,4\,,6\,,8\;,
\label{w1-3-5-oscillator1.2}
\end{equation}
that is $\, w_2=\pi/2\;,\; w_4=3\pi/8\;,\; w_6=5\pi/16\;,\;
w_8=35\pi/128\,,$
the last one
being needed in $\S$~\ref{sec-oscillator3}.

\begin{figure}[h]
\begin{center}
\includegraphics[width=220pt]{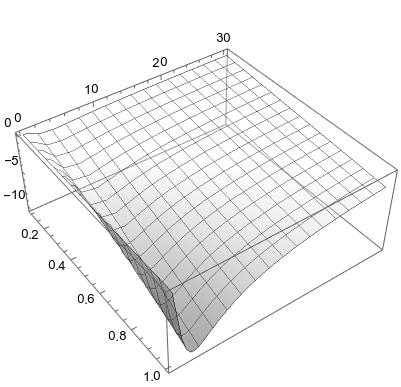}
\includegraphics[width=220pt]{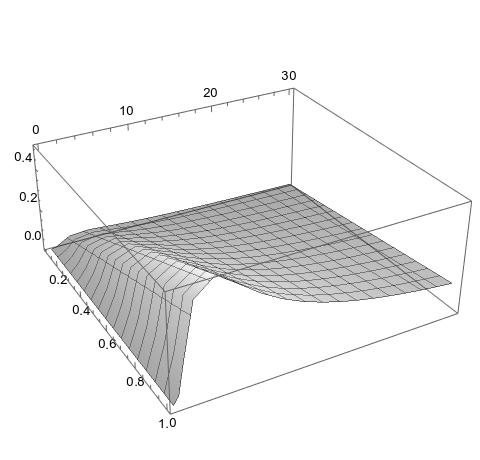}
\caption{Discriminant $\Delta  = 
\Delta(a\,,b) \le 0 \, $ (left) and
coefficient
$\, c_5 =  
c_5(a\,,b) >0 \,$ (right) 
for the quinticated 
form of oscillator
(\ref{v)})
for $\, (a\,,b)\in  ]0\,, 30]\times ]0\,,1]\,$ 
} 
\label{fig-delta-c5-oscillatore2}
\end{center}       
\end{figure}

\noindent
Figure 
\ref{fig-delta-c5-oscillatore2}
depicts
discriminant $\, \Delta = \Delta(a\,,b)\,$ and coefficient 
$\, c_5=c_5(a\,,b)\, $
for varying $\, a>0\, $ and $ \, 0< b\le 1\,,$
showing that condition
(i) of 
\eqref{eqn:conditions} 
is verified;
therefore, 
Theorem \ref{ddue} applies.

\begin{figure}[h]
\begin{center}
\includegraphics[width=220pt]{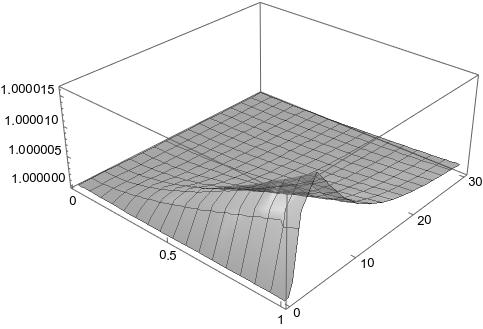}
\includegraphics[width=220pt]{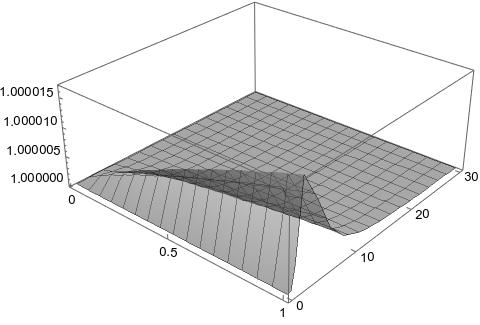}
\caption{Period ratio 
in the case of 
the quinticated form of
oscillator 
\eqref{v)} (left) and oscillator \eqref{vi)} (right) } 
\label{fig-ratio-oscillatore2}
\end{center}       
\end{figure}

\noindent
The qualitative and quantitative behaviour of solution
and period 
for the quinticated 
approximant 
to
oscillator (\ref{v)})
is analogous to that seen in 
$\, \S$~\ref{sec:application2relativistic-oscillator}
for the quinticated relativistic oscillator.
The period ratio stays close to $1\,,$ 
as shown in Figure~\ref{fig-ratio-oscillatore2} (left).
As for the differential operator
\eqref{operator}, 
results were
obtained with parameters $\, a >0\, $ and 
$\, 0< b\le1\,$
and are not reported here, since they
show 
precisely
the
same behaviour and equal order of magnitude
as those
synthetised in 
Figure \ref{fig-norm123}
and Table\ref{tab:norm}
for the case of the relativistic oscillator.

\subsection{Duffing relativistic oscillator (\ref{vi)}) }
\label{sec-oscillator3}

\noindent
The normalised 
IVP, in this case, is:
\begin{equation}
\begin{cases}
\ddot{u}=-u-a^2\;u^3-\dfrac{b\;u}{\sqrt{1+a^2 \;\; u^2}}\;,
\\
u(0)=1\,,\quad \dot{u}(0)=0\;,
\end{cases}
\label{eqn-oscillator3normalised}
\end{equation}
thus 
function
$\, \Phi(u) \,, $ 
defined in   
\eqref{eqn:phi},
becomes:
\begin{equation}
\Phi(u)= \dfrac{1}{2}\; 
(1-u^2)\;(2+a^2+a^2\;u^2)
+
\dfrac{2\,b}{a^2}\, \Big(\sqrt{1+a^2} - \sqrt{1+a^2\; u^2}\Big)\,,
\label{phi-psi-oscillatore3}
\end{equation}
so that forming
$\, \Psi(u)\,$ as in \eqref{main:gensol}
requires
$\,4\, b > 
(2+a^2+a^2\, u^2)\; 
(\sqrt{1+a^2}-\sqrt{1+a^2\,u^2})>0\,,$
being $\, a>0\, $ and $\, 0<u<1\,.$

\vskip 2pt
\noindent
Setting $\, \mathcal{C}= -2/\pi\;$ 
and exploiting formulae 
(\ref{geii})
and
(\ref{w1-3-5-oscillator1.2}),
we have:
\begin{equation}
\begin{aligned}
\alpha_1
&=
\phantom{111}
\mathcal{C}\;  
\left(w_2 + a^2\; w_4 + b\; J_2(a)\right)\;, 
\\
\alpha_3
&= 
- 3\, \mathcal{C}\, \Big( w_2+a^2\, w_4+b\, J_2(a) \Big) 
+ 4\, \mathcal{C}\, \Big( w_4+a^2\, w_6+b\, J_4(a) \Big) 
\\
\alpha_5
&= 
\phantom{-}
5\, \mathcal{C}\, \Big( w_2+a^2\, w_4+b\, J_2(a) \Big) 
- 20\, \mathcal{C}\, \Big( w_4+a^2\, w_6+b\, J_4(a) \Big) 
+ 16\, \mathcal{C}\, \Big( w_6+a^2\, w_8+b\, J_6(a) \Big)\,,
\end{aligned}
\label{c1-3-5-oscillator1.3}
\end{equation}
on which 
$c_1 = -( \alpha_1-3\,\alpha_3+5\, \alpha_5)\;,\;
c_3 = -4 \; ( \alpha_3-5\, \alpha_5)\; $ and
$
c_5 = -16 \;  \alpha_5\;
$
are build, as usual.

\begin{figure}[h]
\begin{center}
\includegraphics[width=250pt]{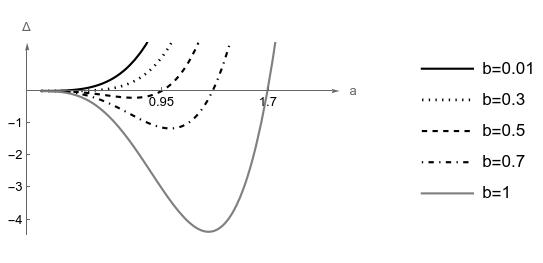}
\caption{Discriminant
$\Delta = \Delta(a\,,b) \, $ 
for the quinticated 
form of oscillator (\ref{vi)})
for various
couples $\,(a\,, b)\,$ } 
\label{fig-delta-oscillatore3}
\end{center}       
\end{figure}

\vskip 2pt
\noindent
Figure 
\ref{fig-delta-oscillatore3}
depicts
discriminant $\, \Delta = \Delta(a\,,b)\,$
for varying $\, a\,, b >0\,, $ 
and
shows how
the values of $\, a \, $ and $\, b\,$ that satisfy either
condition
(i)
or (ii)
of 
\eqref{eqn:conditions} 
are 
closely related.
As an example, 
it is
$\, \Delta(a\,,0.5) \le 0 \,$ and $\, c_5(a\,,0.5)>0\, $ for $\, a \lessapprox 0.95\,;$
similarly,
$\, \Delta(a\,,1) \le 0 \,$ and $\, c_5(a\,,1)>0\, $ for $\, a \lessapprox 1.7\,;$
therefore, in both of these cases,
Theorem \ref{ddue} applies.
In particular,
condition
(i)
starts being significantly verified
when
$\,   b \ge 0.4 \,,$ 
where 
$\,   b \approx 0.4 \,$ 
requires
$a \lessapprox 0.7\,.$
Conversely, for 
$\,   b \lessapprox 0.4\,$ 
and any $\,a\,,$ 
condition
(ii)
is verified and
Theorem \ref{unnno} comes into play.

\begin{figure}[h]
\begin{center}
\includegraphics[width=250pt]{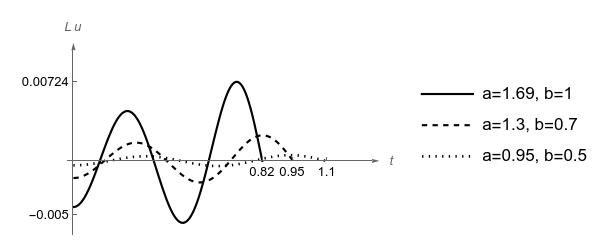}
\caption{$\, L\, u\,,$ 
where $\, u\,$ is solution 
\eqref{udue} with coefficients 
build on 
\eqref{c1-3-5-oscillator1.3},
for the quinticated form of 
oscillator (\ref{vi)}) }
\label{fig-norm-oscillator3}
\end{center}      
\end{figure}

\begin{figure}[h]
\begin{center}
\includegraphics[width=200pt]{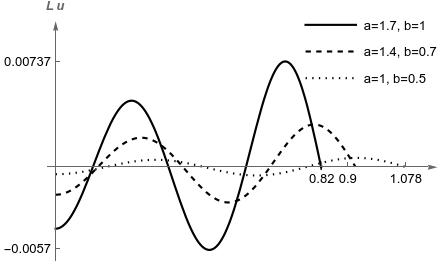}
\caption{$\, L\, u\,,$ 
where $\, u\,$ is solution 
\eqref{uuno} with coefficients 
build on 
\eqref{c1-3-5-oscillator1.3},
for the quinticated form of 
oscillator (\ref{vi)}) }
\label{fig-norm-oscillator3pos}
\end{center}      
\end{figure}

\noindent
The qualitative and quantitative behaviour of solution
and period 
for the quinticated 
approximant to
oscillator (\ref{vi)})
is analogous as those commented
in
$\, \S$~\ref{sec:application2relativistic-oscillator}
and 
$\, \S$~\ref{sec-oscillator2}.
The period ratio stays close to $1\,,$ 
as shown in 
Figure~\ref{fig-ratio-oscillatore2} (right).
As for the differential operator
\eqref{operator}, 
Figure \ref{fig-norm-oscillator3}
and Table\ref{tab:norm-oscillator3}
report results achieved with couples 
$\,(a\,,b)\,$ 
requiring 
the solution $\,u\,$ given by
Theorem
\ref{ddue},
while
Figure \ref{fig-norm-oscillator3pos}
and Table\ref{tab:norm-oscillator3pos}
present the outcome related to 
couples
$\,(a\,,b)\,$ 
for which 
the $\,u\,$
defined in
Theorem 
\ref{unnno}
must be used;
in both cases, 
we attain the same behaviour and equal, or improved, order of magnitude
as 
in the previously studied 
quinticated forms
of oscillators  
\eqref{iii)} and
\eqref{v)}. 

\begin{table}[h]
\centering
\begin{tabular}{ccccccc}
\hline
{a }
& 
{0.95 }
& 
{1.3 }
&
{1.69 }
\\
\hline
{b }
& 
{0.5 }
& 
{0.7 }
&
{1 }
\\
\hline
{$\| L\, u\|_{\infty}$ }
& $\;$ 
{
0.000487249
}
& 
{0.00229373}
& 
{0.00724625}
\\
\hline
\end{tabular}
\caption{Uniform norm of $\, L\, u\,,$ 
with $\, u\,$  solution 
\eqref{udue} and coefficients 
build on 
\eqref{c1-3-5-oscillator1.3},
for the quinticated form of 
oscillator (\ref{vi)}) }
\label{tab:norm-oscillator3}
\end{table}

\begin{table}[h]
\centering
\begin{tabular}{ccccccc}
\hline
{a }
& 
{1 }
& 
{1.4 }
&
{1.7 }
\\
\hline
{b }
& 
{0.5 }
& 
{0.7 }
&
{1 }
\\
\hline
{$\| L\, u\|_{\infty}$ }
& $\;$ 
{
0.00064411
}
& 
{
0.00298016
}
& 
{
0.00737777
}
\\
\hline
\end{tabular}
\caption{Uniform norm of $\, L\, u\,,$ 
with $\, u\,$ solution 
\eqref{uuno}
and coefficients 
build on 
\eqref{c1-3-5-oscillator1.3},
for the quinticated form of 
oscillator (\ref{vi)}) }
\label{tab:norm-oscillator3pos}
\end{table}

\section{Conclusions}
\label{sec:conclusion}

\noindent
In this work
we exploite
\v{C}eby\v{s}\"ev's
fifth-order approximations 
by applying them 
to three popular nonlinear oscillator models,
which share the fact that the integrals obtained 
in the projection are expressible in closed form
by means of complete elliptic integrals of the first and second kind. 
The approximate systems obtained, 
which by their nature constitute very good approximations of the considered models,
are in turn explicitly solved in terms of Jacobian elliptic functions, 
which describe their cosine behaviour. 
The quality of the approximations obtained is confirmed 
in terms of the norm of the deviation of the solution,
and in terms of the ratio between 
the periods
of
the approximating systems 
and the periods 
of the non-approximate systems solutions, 
which are however, in two cases out of three, 
expressible via complete elliptic integrals.
All simulations are performed within the 
{\sl Mathematica} 
scientific environment.

\section*{Acknowledgments}
\noindent
The Authors wish to thank Dr. Mark Sofroniou for many useful discussions.

\section*{Financial disclosure}
\noindent
This research received no external funding.

\section*{Author contributions}
\noindent
The Authors share the content of this work, which is unpublished and has
not been submitted to other journals.
All Authors contributed equally to this work.

\section*{Conflict of interest}
\noindent
The Authors declare no conflict of interests.


\begin{thebibliography}{10}
\expandafter\ifx\csname url\endcsname\relax
  \def\url#1{\texttt{#1}}\fi
\expandafter\ifx\csname urlprefix\endcsname\relax\def\urlprefix{URL }\fi
\expandafter\ifx\csname href\endcsname\relax
  \def\href#1#2{#2} \def\path#1{#1}\fi

\bibitem{maccoll1957}
L.~MacColl, Theory of the relativistic oscillator, American J Physics 25~(8)
  (1957) 535--538.
\newblock \href {https://doi.org/10.1119/1.1934543}
  {\path{doi:10.1119/1.1934543}}.

\bibitem{mickens1998}
R.~Mickens, Periodic solutions of the relativistic harmonic oscillator, J Sound
  and Vibration 212~(5) (1998) 905--908.

\bibitem{younesian2012}
D.~Younesian, H.~Askari, Z.~Saadatnia, M.~KalamiYazdi, Analytical approximate
  solutions for the generalized nonlinear oscillator, Applicable Analysis
  91~(5) (2012) 965--977.

\bibitem{sun2007}
W.~Sun, B.~Wu, C.~Lim, Approximate analytical solutions for oscillation of a
  mass attached to a stretched elastic wire, J Sound and vibration 300~(3--5)
  (2007) 1042--1047.

\bibitem{mickens2001}
R.~Mickens, Mathematical and numerical study of the duffing--harmonic
  oscillator, J Sound and Vibration 244~(3) (2001) 563--567.

\bibitem{van2019}
D.~Van~Hieu, A new approximate solution for a generalized nonlinear oscillator,
  International J Applied and Computational Mathematics 5~(5) (2019) 1--13.

\bibitem{razzak2016analytical}
M.~Razzak, An analytical approximate technique for solving cubic--quintic
  duffing oscillator, Alexandria Engineering J 55~(3) (2016) 2959--2965.

\bibitem{belen2017}
A.~Bel{\'e}ndez, E.~Arribas, T.~Bel{\'e}ndez, C.~Pascual, E.~Gimeno,
  M.~{\'A}lvarez, Closed--form exact solutions for the unforced quintic
  nonlinear oscillator, Advances in Mathematical Physics 2017 (2017).

\bibitem{belen2016}
A.~Bel{\'e}ndez, T.~Bel{\'e}ndez, F.~Martinez, C.~Pascual, M.~Alvarez,
  E.~Arribas, Exact solution for the unforced duffing oscillator with cubic and
  quintic nonlinearities, Nonlinear Dynamics 86~(3) (2016) 1687--1700.

\bibitem{citterio2009}
M.~Citterio, R.~Talamo, The elliptic core of nonlinear oscillators, Meccanica
  44~(6) (2009) 653.

\bibitem{zuniga2013}
A.~Elias-Zuniga, Exact solution of the cubic--quintic duffing oscillator,
  Applied Mathematical Modelling 37~(4) (2013) 2574--2579.

\bibitem{ritelli2010}
G.~Mingari~Scarpello, D.~Ritelli, Exact solution to a first--fifth power
  nonlinear unforced oscillator, Applied Mathematical Sciences 4~(69--72)
  (2010) 3589--3594.

\bibitem{jonckheere1971}
R.~Jonckheere, Determination of the period of nonlinear oscillations by means
  of chebyshev polynomials, Zeitschrift fur angewandte Mathematik und Mechanik
  51~(5) (1971) 389--393.

\bibitem{zuniga2014}
A.~Elias-Zuniga, Quintication method to obtain approximate analytical solutions
  of non--linear oscillators, Applied Mathematics and Computation 243 (2014)
  849--855.

\bibitem{biglabo2019}
A.~Big-Alabo, Approximate period for large--amplitude oscillations of a simple
  pendulum based on quintication of the restoring force, European Journal of
  Physics 41~(1) (2019) 015001.

\bibitem{belen2009}
A.~Bel{\'e}ndez, D.~M{\'e}ndez, M.~Alvarez, C.~Pascual, T.~Bel{\'e}ndez,
  Approximate analytical solutions for the relativistic oscillator using a
  linearized harmonic balance method, International J Modern Physics B 23~(04)
  (2009) 521--536.

\bibitem{biazar2014}
J.~Biazar, M.~Hosami, An easy trick to a periodic solution of relativistic
  harmonic oscillator, J Egyptian Mathematical Society 22~(1) (2014) 45--49.

\bibitem{hosen2016}
M.~Hosen, M.~Chowdhury, M.~Ali, A.~Ismail, A novel analytical approximation
  technique for highly nonlinear oscillators based on the energy balance
  method, Results in Physics 6 (2016) 496--504.

\bibitem{nayfeh1973}
A.~Nayfeh, Perturbation Methods, Wiley \& Sons, New York, USA, 1973.

\bibitem{nayfeh1979}
A.~Nayfeh, D.~Mook, Nonlinear oscillations, John Wiley \& Sons, New York, USA,
  1979.

\bibitem{mickens1996}
R.~Mickens, Oscillations in planar dynamic systems, World Scientific,
  Singapore, 1996.

\bibitem{krylov1949}
N.~Krylov, N.~Bogoliubov, Introduction to non--linear mechanics, Princeton
  University Press, Princeton, NJ, USA, 1949.

\bibitem{he2002}
J.~He, Preliminary report on the energy balance for nonlinear oscillations,
  Mechanics Research Communications 29~(2--3) (2002) 107--111.

\bibitem{belendez2008a}
A.~Bel{\'e}ndez, T.~Bel{\'e}ndez, A.~M{\'a}rquez, C.~Neipp, Application of
  {H}e’s homotopy perturbation method to conservative truly nonlinear
  oscillators, Chaos, Solitons \& Fractals 37~(3) (2008) 770--780.

\bibitem{mickens2010}
R.~Mickens, Truly nonlinear oscillations: harmonic balance, parameter
  expansions, iteration, and averaging methods, World Scientific, Singapore,
  2010.

\bibitem{gottlieb2004}
H.~Gottlieb, Harmonic balance approach to periodic solutions of non--linear
  jerk equations, J Sound and Vibration 271~(3--5) (2004) 671--683.

\bibitem{gottlieb2006}
H.~Gottlieb, Harmonic balance approach to limit cycles for nonlinear jerk
  equations, J Sound and Vibration 297~(1--2) (2006) 243--250.

\bibitem{wu2006}
B.~Wu, C.~Lim, W.~Sun, Improved harmonic balance approach to periodic solutions
  of non--linear jerk equations, Physics Letters A 354~(1--2) (2006) 95--100.

\bibitem{belen2008}
A.~Bel{\'e}ndez, D.~M{\'e}ndez, T.~Bel{\'e}ndez, A.~Hern{\'a}ndez, M.~Alvarez,
  Harmonic balance approaches to the nonlinear oscillators in which the
  restoring force is inversely proportional to the dependent variable, J Sound
  and Vibration 314~(3--5) (2008) 775--782.

\bibitem{kovacic2011}
I.~Kovacic, M.~Brennan, The Duffing equation: nonlinear oscillators and their
  behaviour, John Wiley \& Sons, New York, USA, 2011.

\bibitem{cveticanin2014}
L.~Cveti\'canin, Strong Nonlinear {O}scillators, Springer, Cham, Switzerland,
  2014.

\bibitem{rothe1985}
F.~Rothe, The periods of the volterra--lotka system, J Reine Angew. Math 355
  (1985) 129--138.

\bibitem{ritelli2004}
S.~Foschi, G.~Mingari~Scarpello, D.~Ritelli, Higher order approximation of the
  period--energy function for single degree of freedom hamiltonian systems,
  Meccanica 39~(4) (2004) 357--368.

\bibitem{azami2009}
R.~Azami, D.~Ganji, H.~Babazadeh, A.~Dvavodi, S.~Ganji, He's max--min method
  for the relativistic oscillator and high order duffing equation,
  International J Modern Physics B 23~(32) (2009) 5915--5927.

\bibitem{belendez2008higher}
A.~Bel{\'e}ndez, C.~Pascual, E.~Fern{\'a}ndez, C.~Neipp, T.~Bel{\'e}ndez,
  Higher--order approximate solutions to the relativistic and duffing--harmonic
  oscillators by modified {H}e's homotopy methods, Physica Scripta 77~(2)
  (2008) 025004.

\bibitem{jamshidi2010application}
N.~Jamshidi, D.~Ganji, Application of energy balance method and variational
  iteration method to an oscillation of a mass attached to a stretched elastic
  wire, Current Applied Physics 10~(2) (2010) 484--486.

\bibitem{zhao2009he}
L.~Zhao, He’s frequency--amplitude formulation for nonlinear oscillators with
  an irrational force, Computers \& Mathematics with Applications 58~(11--12)
  (2009) 2477--2479.

\bibitem{belendez2007application}
A.~Bel{\'e}ndez, A.~Hern{\'a}ndez, T.~Bel{\'e}ndez, M.~Alvarez, S.~Gallego,
  M.~Ortuno, C.~Neipp, Application of the harmonic balance method to a
  nonlinear oscillator typified by a mass attached to a stretched wire, J Sound
  and Vibration 302~(4--5) (2007) 1018--1029.

\bibitem{marion2013classical}
J.~Marion, Classical dynamics of particles and systems, Academic Press, New
  York, USA, 2013.

\bibitem{graham-knuth-patashnik}
R.~Graham, D.~Knuth, O.~Parashnik, Concrete Mathematics, 2nd ed.,
  Addison--Wesley, Reading, MASS, USA, 1994.

\bibitem{trefethen}
L.~Trefethen, Approximation Theory and Approximation Practice, SIAM,
  Philadelphia, PA, USA, 2019.

\bibitem{phillips-taylor}
G.~Phillips, P.~Taylor, Theory and Applications of Numerical Analysis, 2nd ed.,
  Academic Press, Elsevier Science $\&$ Technology, Boston, MASS, USA, 1996.

\bibitem{belen2012}
A.~Bel{\'e}ndez, M.~Alvarez, J.~Franc{\'e}s, S.~Bleda, T.~Bel{\'e}ndez,
  A.~N{\'a}jera, E.~Arribas, Analytical approximate solutions for the
  cubic--quintic duffing oscillator in terms of elementary functions, J Applied
  Mathematics Volume 2012, Article ID 286290 (2012).

\bibitem{belen2016b}
A.~Bel{\'e}ndez, A.~Hernandez, T.~Bel{\'e}ndez, C.~Pascual, , M.~Alvarez,
  E.~Arribas, Solutions for conservative nonlinear oscillators using an
  approximate method based on chebyshev series expansion of the restoring
  force, ACTA PHYSICA POLONICA A 130~(3) (2016) 667--678.

\bibitem{khalil2021}
H.~Khalil, M.~Khalil, I.~Hashim, P.~Agarwal, Extension of operational matrix
  technique for the solution of nonlinear system of caputo fractional
  differential equations subjected to integral type boundary constrains,
  Entropy 29~(3) (2021) 1154.

\bibitem{grad2000}
I.~Gradshteyn, J.~Ryzhik, Table of Integrals, {S}eries and {P}roducts 6th ed,
  Academic Press, New York, USA, 2000.

\bibitem{byrd1971}
P.~Byrd, M.~Friedman, Handbook of elliptic integrals for engineers and
  scientists, Springer Berlin, New York, USA, 1971.

\bibitem{wri}
S.~Wolfram, An Elementary Introduction to the Wolfram Language, 2nd ed.,
  Wolfram Media, Inc., Urbana--Champaign, ILL, USA, 2017.

\end{thebibliography}

\end{document}